# Parameterization adaption for 3D shape optimization in aerodynamics

Badr Abou El Majd
Computer Science and Decision Aiding Laboratory
Hassan II University
B.P. 5366 Maârif, Casablanca, Morocco
b.abouelmajd@fsac.ac.ma

**Abstract**

When solving a PDE problem numerically, a certain mesh-refinement process is always implicit, and very classically, mesh adaptivity is a very effective means to accelerate grid convergence. Similarly, when optimizing a shape by means of an explicit geometrical representation, it is natural to seek for an analogous concept of parameterization adaptivity. We propose here an adaptive parameterization for three-dimensional optimum design in aerodynamics by using the so-called "Free-Form Deformation" approach based on 3D tensorial Bézier parameterization. The proposed procedure leads to efficient numerical simulations with highly reduced computational costs.

## 1. Introduction

Within few years, numerical shape optimization is playing a great role in aerodynamic aircraft design. It enables to design and improve the shape of some or all of the components of the aircraft by minimising a cost functional subject to physical and geometrical constraints. This cost function trusts in the prior solution of a complex set of partial-differential equations (PDEs), such as those governing compressible aerodynamics (e.g. the Euler equations). Whence, the optimization process suffers from the high computational effort for the flow simulations around 3D configurations when the accuracy requirement is high. Thus, our efforts is mostly concentrated on improving the convergence rate of the numerical procedures both from the viewpoint of cost-efficiency and accuracy by handling the parametrization of the shape to be optimized.

When solving a PDE problem numerically, a certain mesh-refinement process is always implicit, and very classically, mesh adaptivity instead of, or in conjunction with increasing the number of degrees of freedom, is a very effective means to accelerate grid convergence. Similarly, when optimizing a shape by means of an explicit geometrical representation, as we advocate, it is natural to seek for an analogous concept of parameterization adaptivity. We propose here a self-adaptive procedure for a three-dimensional optimum-design in aerodynamics by using the so-called Free-Form Deformation (FFD) method [12]. This approach is studied initially in the framework of the Bézier parametrization and applied to a geometrical arc reconstruction [3].

This paper is organized as follows. sections 2.1 introduces some properties of the classical Bézier paramerization. Then, we recall the concept of Free-Form Deformation which allows to extend the concept of shape representation to three-dimensional cases. In section 3, we present the notion of parametrization adaption within the framework of FFD approach. We apply the self-adaptivity approach to optimum shape design in 3D aerodynamics. Finally, we conclude and give some prespectives.

# 2. Shape representation

## 2.1 Bézier parameterization

We begin with the simplest situation of a two-dimensional geometry for which we employ a Bézier shape representation:

$$x(t) = \sum_{k=0}^{n} B_n^k(t)\, x_k\,, \quad y(t) = \sum_{k=0}^{n} B_n^k(t)\, y_k$$

in which the parameter $t$ varies from 0 to 1, $n$ is the degree of the parameterization,

$$B_n^k(t) = C_n^k\, t^k\, (1-t)^{n-k}$$

is a Bernstein polynomial, $C_k^n = \frac{n!}{k!(n-k)!}$, and

$$P_k = \binom{x_k}{y_k} \qquad (k = 0,1,\ldots,n)$$

is the generic control point. The coordinates of these control points are split into two vectors

$$X = \{x_k\}\,, \quad Y = \{y_k\}\,, \quad k = 0,1,\ldots,n,$$

and we refer to the vector $X$ as the *support* of the parameterization, and the vector $Y$ as the *design* vector. Typically, we optimize the design vector for fixed support according to some physical criterion, such as drag reduction in aerodynamics. The somewhat unsymmetrical roles dispensed to the vectors $X$ and $Y$ are chosen to reduce (to $n$ essentially) the dimension of the search space in the optimization phase, which is the most numerically costly and subject to numerical stiffness.
We also use the notation:

$$x(t) = B_n(t)^T\, X\,, \quad y(t) = B_n(t)^T\, Y\,,$$

in which the vector $B_n(t)^T = (B_n^0(t), B_n^1(t), \ldots, B_n^n(t))$. In all this article, only supports for which the sequence $\{x_k\}$ is monotone increasing are said to be admissible and considered throughout. Thus, the function $x(t)$ is monotone-increasing and defines a one-to-one mapping of, say, [0,1] onto itself.

## 2.2 Free-Form Deformation approach

A critical issue in aerodynamic design is the choice of the shape parameterization. Parameterization techniques for practical 3D aerodynamic shape optimization have to fulfill several criteria :

- the parameterization should be able to take into account complex geometries, possibly including constraints and singularities ;
- the number of parameters should be as small as possible, since the stiffness of the shape optimization numerical formulation increases abruptly with the number of parameters;
- the parameterization should allow to control the smoothness of the resulting shapes.

A survey of shape parameterization techniques for multi-disciplinary optimization, which are analyzed according to the previous criteria, is proposed by Samareh [6]. Following his recommendation, conclusions, the Free-Form Deformation (FFD) [12] technique is adopted in the present study, since it provides an easy and powerful framework for the deformation of complex shapes, such as generic or elaborate aerodynamic configuration.

The FFD technique originates from the Computer Graphics field [12]. It allows the deformation of an object in a 2D or 3D space, regardless of the representation of this object. Instead of manipulating the surface of the object directly, by using classical B-Splines or Bézier parameterization of the surface, the FFD technique defines a deformation field over the space embedded in a lattice which is built around the object. By modifying the space coordinates inside the lattice, the FFD technique deforms the object, regardless of its geometrical description. In particular, the initial geometry, in our applications, is usually defined by a general, Finite-Element-type unstructured simplicial grid.

More precisely, consider a three-dimensional hexaedral lattice embedding the object to be deformed. Figure [fig:ffd1] shows an example of such a lattice built around a typical wing. A local coordinate system $(\xi, \eta, \zeta)$ is defined in the lattice, with $(\xi, \eta, \zeta) \in [0,1] \times [0,1] \times [0,1]$. As a result of the deformation, the displacement $\Delta q$ of each point $q$ inside the lattice is here defined by a third-order Bézier tensor product:

$$\Delta q = \sum_{i=0}^{n_i} \sum_{j=0}^{n_j} \sum_{k=0}^{n_k} B_{n_i}^{i}(\xi_q)\, B_{n_j}^{j}(\eta_q)\, B_{n_k}^{k}(\zeta_q)\, \Delta P_{ijk}.$$

$B_{n_i}^{i}$, $B_{n_j}^{j}$ and $B_{n_k}^{k}$ are again Bernstein polynomials of order $n_i$, $n_j$ and $n_k$. $(\Delta P_{ijk})_{0 \le i \le n_i, 0 \le j \le n_j, 0 \le k \le n_k}$ are weighting coefficients, or control points displacements, which are used to monitor the deformation and are considered as design variables during the shape optimization procedure. The critical point is that only the *shape deformation* is represented not the shape itself.

This technique is illustrated by Figure 1. A lattice is built around a wing and a Bézier tensor product of degree$n_i = 4$, $n_j = 1$ and $n_k = 1$ is defined over this lattice. Corner control points (filled markers) are supposed to be frozen in order to keep leading and trailing edges fixed during the deformation, whereas other control points (empty markers) are allowed to move vertically (Figure 1(a)). When these control points are moved, their displacements define a continuous deformation inside the lattice according to [6], yielding a shape deformation. The deformed lattice and shape can be seen in Figure 1(b).

The FFD technique described above is well suited to complex shape optimization, thanks to the following properties:

- the initial shape can be exactly represented (no deformation occurs when all weighting coefficients are zero);
- the deformation is performed whatever the complexity of the shape (this is a free-form technique);
- geometric singularities can be taken into account (the initial shape including its singularities is deformed);
- the smoothness of the deformation is controlled (the deformation is ruled by Bernstein polynomials);
- the number of design variables depends on the user's choice (the deformation is independent of the shape itself);

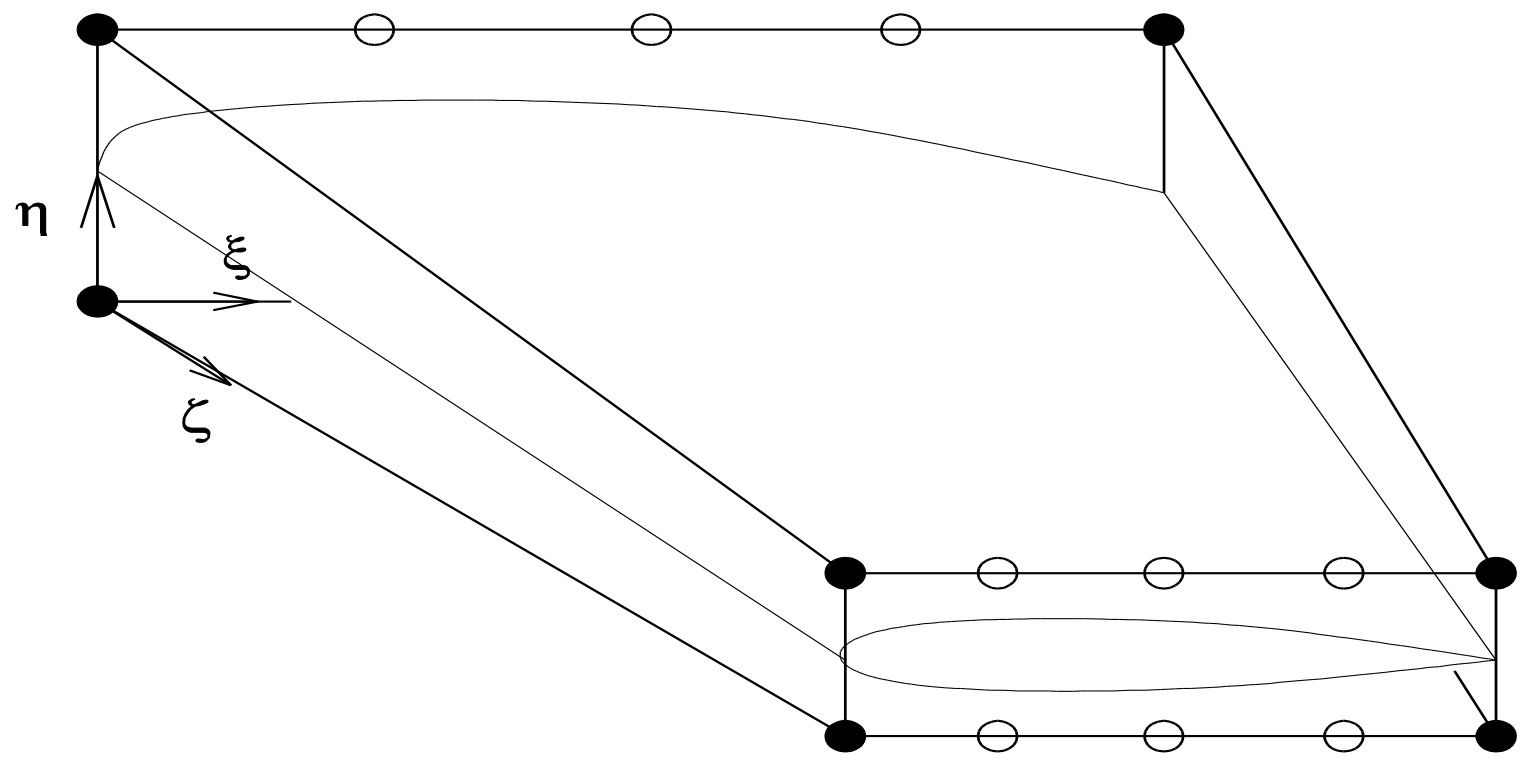

(a) initial FFD lattice 4-1-1

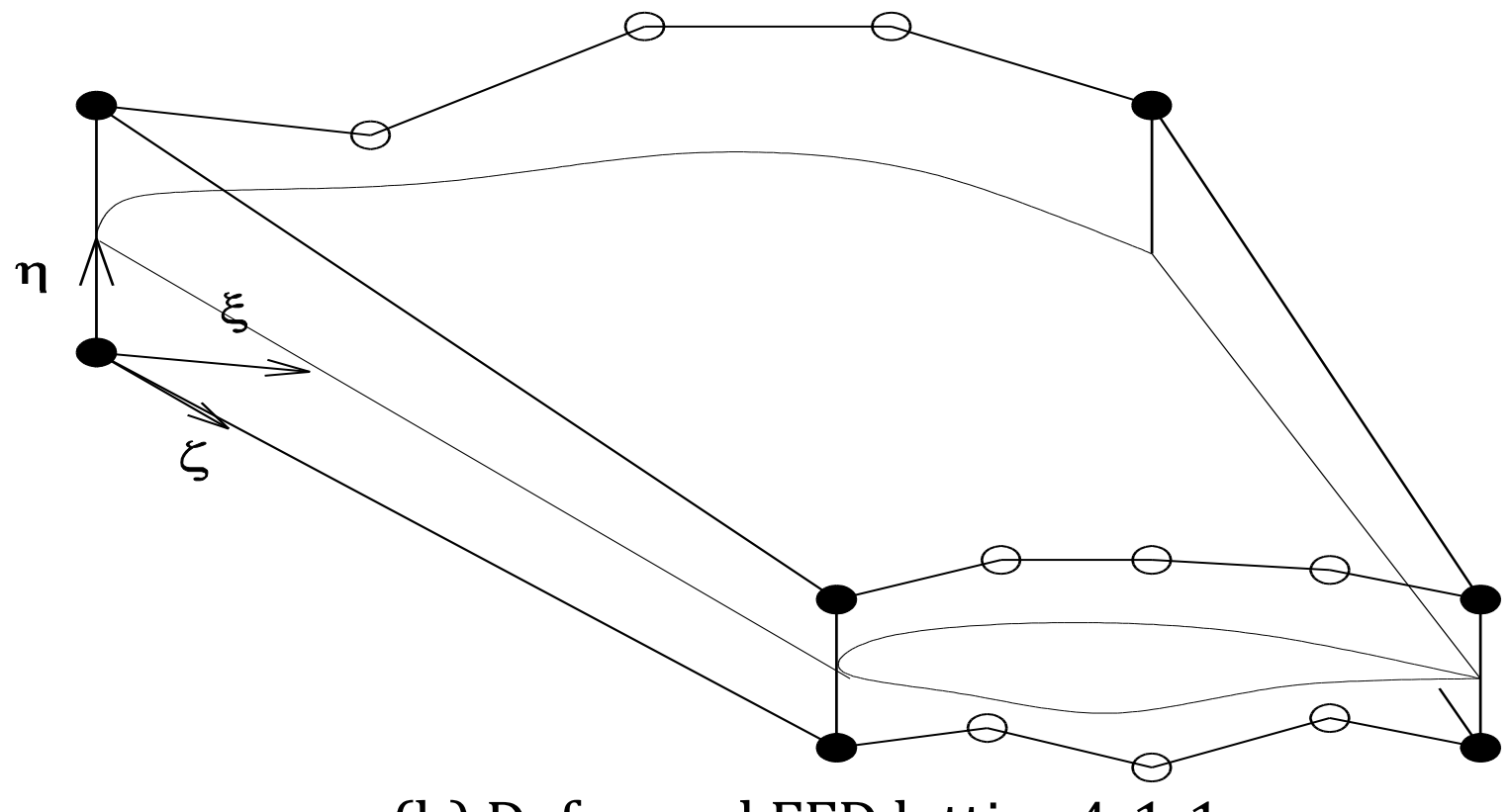

(b) Deformed FFD lattice 4-1-1

Figure 1: Example of Free-Form Deformation: by moving some control points of the lattice, a deformation field is defined continuously inside the lattice, yielding a shape deformation.

# 3. 3D parameterization adaptivity

## 3.1 Motivations

Before defining our concept of parameterization adaptivity, we discuss some elements that has motivated its construction. For this purpose, we use an intrinsic formulation of shape-reconstruction problem, initially introduced in [5][11][12]:

$$\min_{\gamma} \Im(\gamma) := \int_{\gamma} \frac{1}{2} [y(x) - \overline{y}(x)]^2 dx$$

where $\gamma$ is the unknown shape analytically represented by $y(x)$; $\overline{y}(x)$ is the analogous analytical representation of a given target curve $\overline{\gamma}(x)$, subsequently assumed, without great loss of generality, to be a Bézier curve of degree $n$ and support $X$. This problem is transformed into a parametric optimization by assuming Bézier representations of the curves over the support $X$:

$$\min_{Y \in \mathbf{R}^{n+1}} j_n(Y) := \int_0^1 \frac{1}{2} \left[B_n(t)^T (Y - \overline{Y})\right]^2 n B_{n-1}(t)^T \Delta X dt$$

The symbol $\Delta$ represents the forward-difference operator that appears when differentiating Bernstein polynomials.

Since the functional is quadratic, the parametric gradient is linear (in $Y$):

$$j'_n(Y) = A(X)Y - b(X)$$

where

$$A(X) = \int_0^1 B_n(t) B_n(t)^T n B_{n-1}(t)^T \Delta X dt$$

and

$$b(X) = \int_0^1 B_n(t) B_n(t)^T \overline{Y} n B_{n-1}(t)^T \Delta X dt$$

In particular, for a uniform support $X$, the matrix $A$ reduces to the simple form :

$$A(X) = \int_0^1 B_n(t) B_n(t)^T dt = Aij$$

in which the coefficients $Aij$ are obtained by a simple calculation :

$$Aij = \frac{1}{2n+1} \frac{C_n^i C_n^j}{C_{2n}^{i+j}}$$

For this shape-inverse problem, the optimization problem is equivalent to solving the linear system ($X$ is fixed during the optimization process),

$$AY = b(X)$$

The matrix $A$ indicates how does the parametrization condition the stifness of the optimization iteration. By computing the condition number for differente values of the parametrization degree $n$, we observe, as shown in figure 2, that the condition number of the matrix $A$ increase with $n$. Because, in practice, the parameter $n$ must be sufficienty fine subject to obtain efficient solution, the linear system $AY = b(X)$ is ill-conditioned due to the cluster of small singular values of the matrix $A$ (as depicted in figure 3). Thus, the vector $\tilde{Y} = A^{-1}b(X)$ ($A^{-1}$ is the inverse of $A$) is a usually meaningless bad approximation to the exact solution $Y$. Hence, the so-called regularization techniques are needed to stabilize such ill-conditioned problem and obtain meaningful solution estimates.

Undoubtedly, the most common and well-known form of regularization is the one known as Tikhonov regularization. The idea is to seek the regularized solution $Y_\rho$ as the minimizer of the following weighted functional

$$\phi_\rho(Y) = \|AY - b\|_2^2 + \rho\|Y\|_2^2,$$

where the first term corresponds to the residual norm, and the second to a side constraint imposed on the solution. The regularization parameter $\rho$ is an important quantity which controls the properties of the regularized solution, and $\rho$ should therefore be chosen with care.

Intuitively, the Tikhonov regularization tries to find a good trade-off between two requirements:

1. $Y_\rho$ should give a small residual $AY_\rho - b$.
2. $Y_\rho$ should be regularized with respect to the 2-norm.

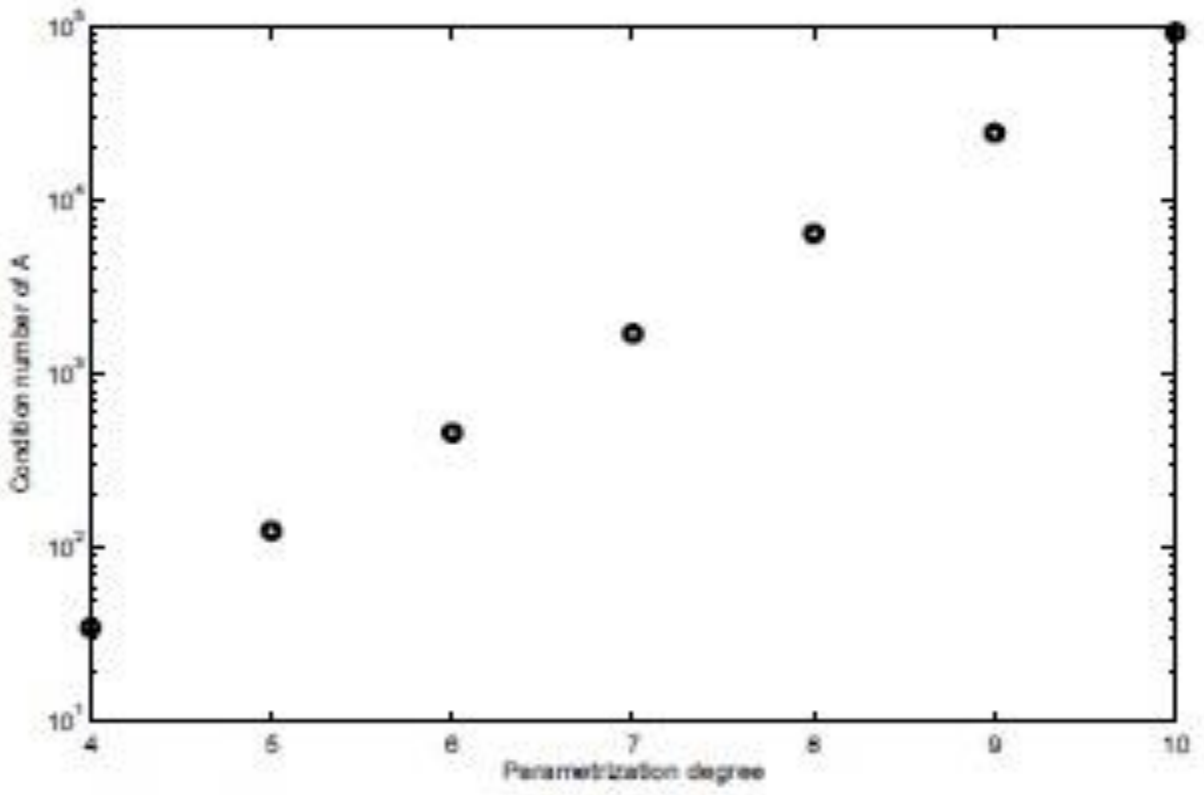

Figure 2: The condition number of A for different values of the parameterization degree

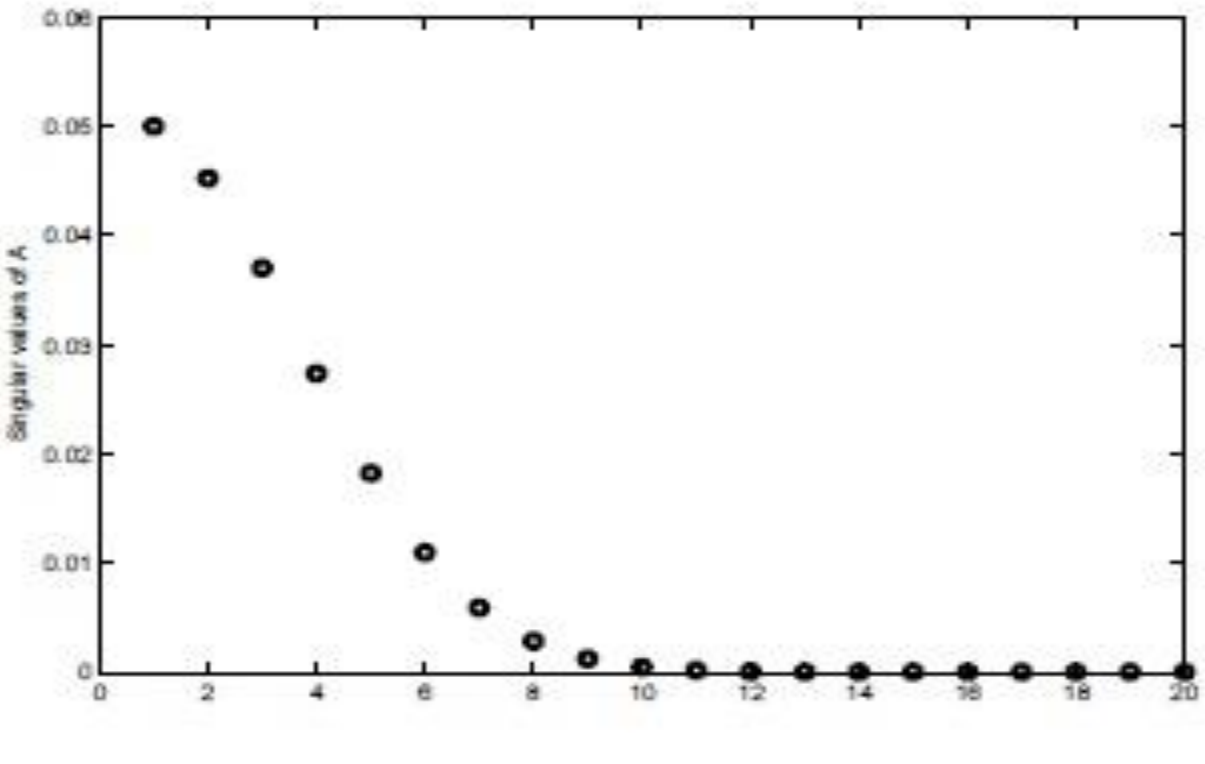

Figure 3: The singular values of the matrix A This above reasoning explains why, for this shape-recontruction problem, the regularization of the parametrization is necessary to reduce the stifness of the optimization problem. In two-dimentional cases, our algorithm require two complementary phases:

1. Optimization: optimize the design vector $Y$ for fixed support $X = X^0$ according to some criterion; let $Y^0$ be the result of this phase.
2. Regularization: given the parametrization $(X^0, Y^0)$ of an approximate optimum shape $\gamma^0$, the new support $X^1$ is taken to be the better support for which the total variation (TV) in the components of the corresponding vector $Y^1$ is minimal, such that the correspondant shape $\gamma^1$ to $(X^1, Y^1)$ approximates $\gamma^0$ in the sense of least squares; subtitute $X^1$ to $X^0$.

We note that the support of the parametrization is the regularization parameter; it plays the same role as the parameter $\rho$ in Tikhonov regularization.

## 3.2 Principles

For complex three-dimentional problems, such as those encountered in aerodynamics, the Free-Form Deformation (FFD) approach is adopted. In this case, one proposes here the notion of parametrization adaptivity. In the framework of aerodynamic shape optimization, two outcomes are expected:

- reach a shape of better fitness (i.e. decrease the distance between reachable shapes and the best existing shape) ;
- increase the convergence rate (i.e. improve the conditioning of the numerical optimization problem by modifying the topology of the cost function).

Representing the shape by Free-Form Deformation allows to deform mesh and shape simultaneously. The quality of the shape is kept and a costly re-meshing is avoided. During the optimization, The shape together with the grid embedded into the FFD control volume are deformed. Continuity with the CFD mesh outside of the volume is guaranteed by fixing the control point of the boundary.

The parameterization adaptivity in this general context is inspired from the approach in two-dimensional presented in the previous section. The main idea consiste to yield the control volume more regular during the shape optimization process. In fact, it was observed that the control volume become very irregular after some iterations. So, we use the parametrization adaptivity process, explained below, in order to change the actual neighboring control points by a new one which is more regularized.

Denote by $(X_{min}, Y_{min}, Z_{min})$ and $(X_{max}, Y_{max}, Z_{max})$ the corners of a deformation region. Let $(x^0, y^0, z^0)$ a Cartesian coordinates of an interior point $Q^0$ of the initial mesh. The deformation of the lattice around an object is specified by changing the initial control points ${P^0}_{ijk}$ defined by,

$$ {P^0}_{ijk} = \left( X_{min} + \frac{i}{n}(X_{max} - X_{min}, Y_{min} + \frac{j}{m}(Y_{max} - Y_{min}, Z_{min} + \frac{k}{l}(Z_{max} - Z_{min}) \right) $$

In order to find the current mesh before fitness evaluation, we add each node $Q^0$ of the initial mesh with the corresponding deformation, it follows that

$$Q = Q^0 + \sum_{i=0}^{n}\sum_{j=0}^{m}\sum_{k=0}^{l} B_n^i\,(\xi_{Q^0})\,B_m^j(\eta_{Q^0})\,B_l^k(\zeta_{Q^0})\,{P^0}_{ijk}$$

where,

$$\xi_{Q^0} = \frac{x^0 - X_{min}}{X_{max} - X_{min}}, \quad \eta_{Q^0} = \frac{y^0 - Y_{min}}{Y_{max} - Y_{min}}, \zeta_{Q^0} = \frac{z^0 - Z_{min}}{Z_{max} - Z_{min}}.$$

After performing $k$ optimization steps, the shape and the computational grid follow all deformation applied to the control volume; the resulting control volume ${P^k}_{ijk}$ is irregular as illustrated by Figure 4. The current mesh inside the lattice volume is achieved by adding the initial mesh to the current control points,

$$Q^k = Q^0 + \sum_{i=0}^{n}\sum_{j=0}^{m}\sum_{k=0}^{l} B_n^i\,(\xi_Q)\,B_m^j(\eta_Q)\,B_l^k(\zeta_Q)\,{P^k}_{ijk}$$

Our approach of parametrization adaptivity consist to restart the optimization by both the current mesh and the initial volume ${P^0}_{ijk}$, it follows that for each node $Q(x, y, z)$ of the current mesh,

$$Q^k = Q^k + \sum_{i=0}^{n}\sum_{j=0}^{m}\sum_{k=0}^{l} B_n^i\,(\xi_Q)\,B_m^j(\eta_Q)\,B_l^k(\zeta_Q)\,{P^0}_{ijk}$$

where,

$$\xi_Q = \frac{x - X_{min}}{X_{max} - X_{min}}, \quad \eta_Q = \frac{y - Y_{min}}{Y_{max} - Y_{min}}, \quad \zeta_Q = \frac{z - Z_{min}}{Z_{max} - Z_{min}}.$$

So, the optimization restarts with exactly the same deformed mesh/shape but a new regular control volume as depicted in Figure 4. This adaption pocedure happens at some steps during the optimization process. It is not clear practically when the adaptivity procedure takes place optimally.

# 4. Results

## 4.1 Test-case description

The test-case considered here corresponds to the optimization of the wing shape of a business aircraft (courtesy of Piaggio Aero Industries) in a transonic regime. The free-stream Mach number is $M_\infty = 0.83$ and the incidence $\alpha = 2°$. Initially, the wing section corresponds to the NACA 0012 airfoil. An unstructured mesh, composed of 31124 nodes and 173 445 elements, is generated around the wing, including a refined area in the vicinity of the shock (Figure 5). Flow fields are obtained by solving compressible Euler equations using a finite-volume method.

The goal of the optimization is to reduce the drag coefficient $C_D$ subject to the constraint that the lift coefficient $C_L$ should not decrease more than 0.1%. The constraint is taken into account using a penalization approach. Then, the resulting cost function is :

$$\mathcal{J}_{OPT} = \frac{C_D}{C_{D_0}} + 10^4 \max(0, 0.999 - \frac{C_L}{C_{L_0}}).$$

$C_{D_0}$ and $C_{L_0}$ are respectively the drag and lift coefficients corresponding to the initial shape (NACA 0012 section).

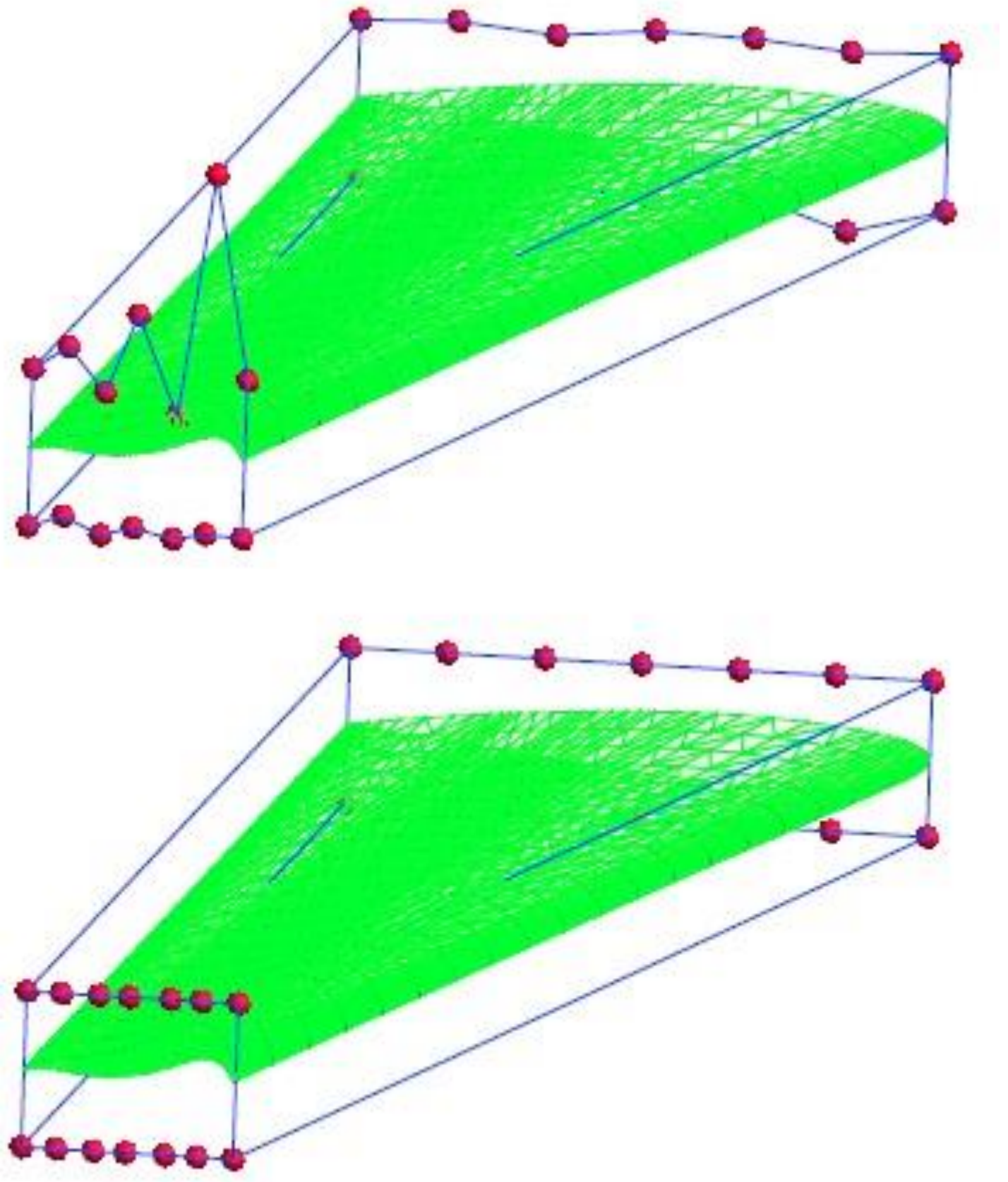

Figure 4: The control points of the lattice before and after adaption.

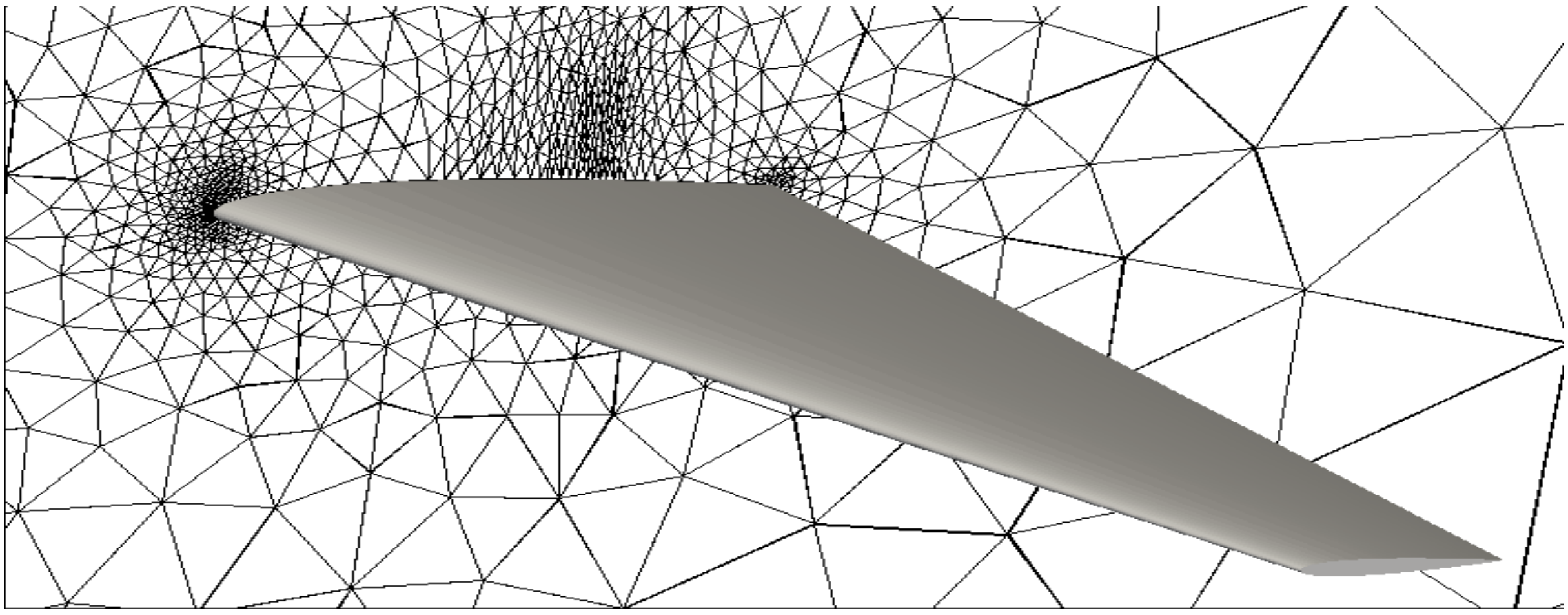

Figure 5: Initial wing shape and mesh in the symmetry plane.

The FFD lattice is built around the wing with $\xi$, $\eta$ and $\zeta$ in the chord-wise, span-wise and thickness directions respectively. The lattice is chosen in order to fit the planform of the wing. Then, the leading and trailing edges are kept fixed during the optimization by freezing the control points that correspond to $i = 0$ and $i = n_i$. Moreover, control points are only moved vertically. Results are presented for three parameterizations. The coarsest one corresponds to $n_i = 3$, $n_j = 1$ and $n_k = 1$. Therefore, $(4 - 2) \times 2 \times 2 = 8$ degrees of freedom are taken into account in the optimization. The medium parameterization corresponds to $n_i = 6$, $n_j = 1$ and $n_k = 1$ and counts $(7 - 2) \times 2 \times 2 = 20$ degrees of freedom. Finally, the finest parameterization corresponds to $n_i = 9$, $n_j = 1$ and $n_k = 1$ and counts $(10 - 2) \times 2 \times 2 = 32$ degrees of freedom. In this study, The Nelder-Mead simplex method [7] is used as optimization algorithm.

For the coarse parameterization $(3 \times 1 \times 1)$ two stratgies are compared:

- Basic method : optimization without adaption until full convergence;
- Adaptive method: optimization with adaption.

In this test case, the adaption process occurs after 100 iterations of the shape optimization process.

## 4.2 Aerodynamic coefficients

The aerodynamic coefficients obtained for each method are compared in Table 1. The lift coefficient is approximately maintained or slightly increased by the shape optimization process. Important reductions of the drag coefficient are reported. As observed, the utilization of adaption method improve significantly the aerodynamic performance.

| Method | $C_L$ | $C_D$ | Cost |
|---|---|---|---|
| Reference | 0.319192893 | 0.026352608 | 1. |
| Basic method | 0.318874966 | 0.017450289 | 0.662184501 |
| Adaptive method | 0.318999078 | 0.016299483 | 0.618515468 |

Table 1: Comparison of aerodynamic coefficients and cost function values.

## 4.3 Convergence history plots

Figure 6 shows a comparison of the convergence for the two strategies under consideration. The adaptive method is significantly efficient than the basic method, yielding a shape of better fitness using a smaller computational effort.

## 4.4 Flows

A comparison of the flow fields for the final shapes obtained with the different strategies is presented in Figures 7 to 9. The Mach number field on the wing surface and Mach number contours in the symmetry plane are represented. Visibly, this drag reduction exercise results in a strong reduction of the shock wave. Using a basic method, the shock reduction is not as important, whereas in the adaption approaches, the shock at the root section disappears.

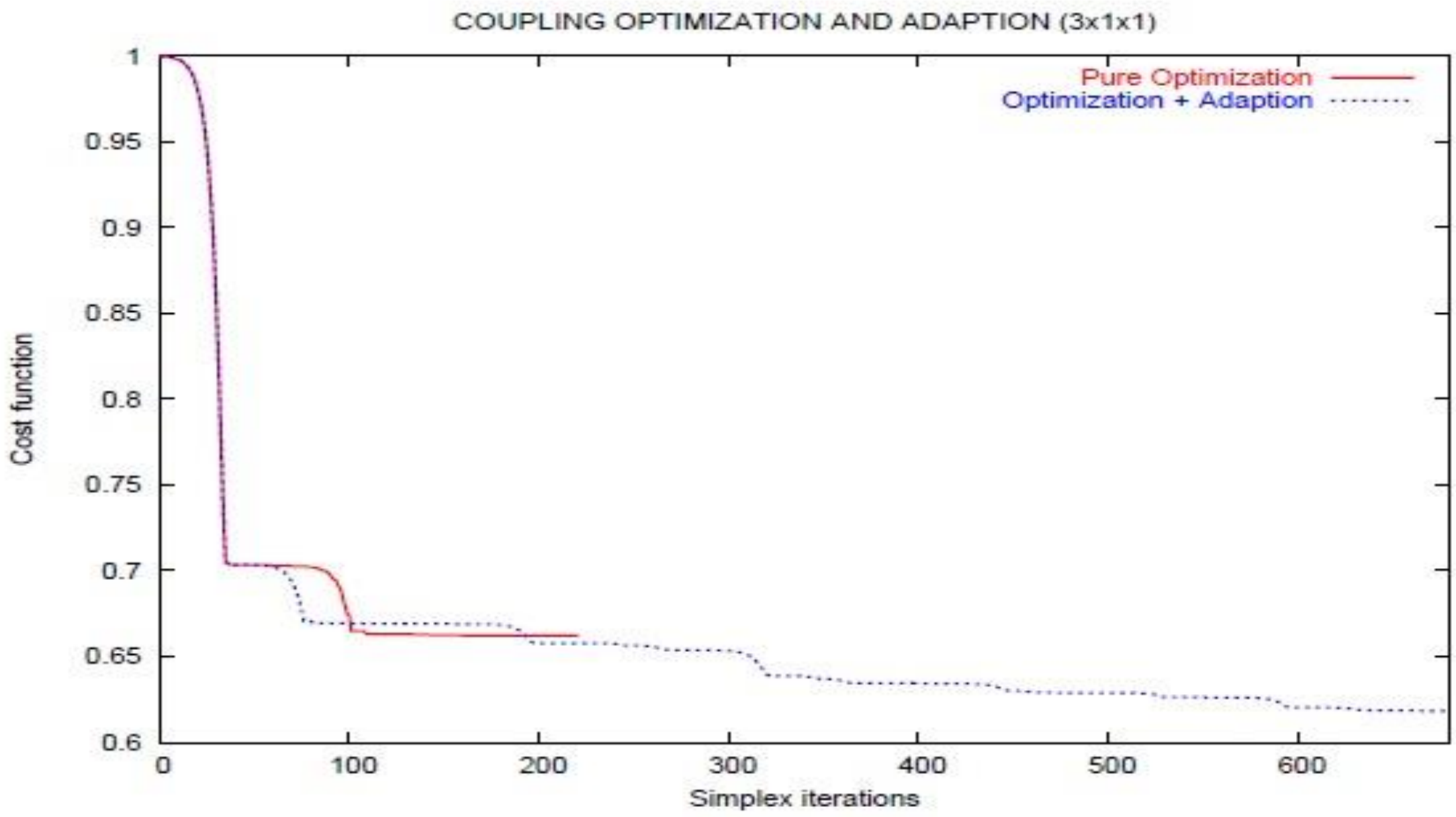

Figure 6: Comparison of the convergence history for the two strategies.

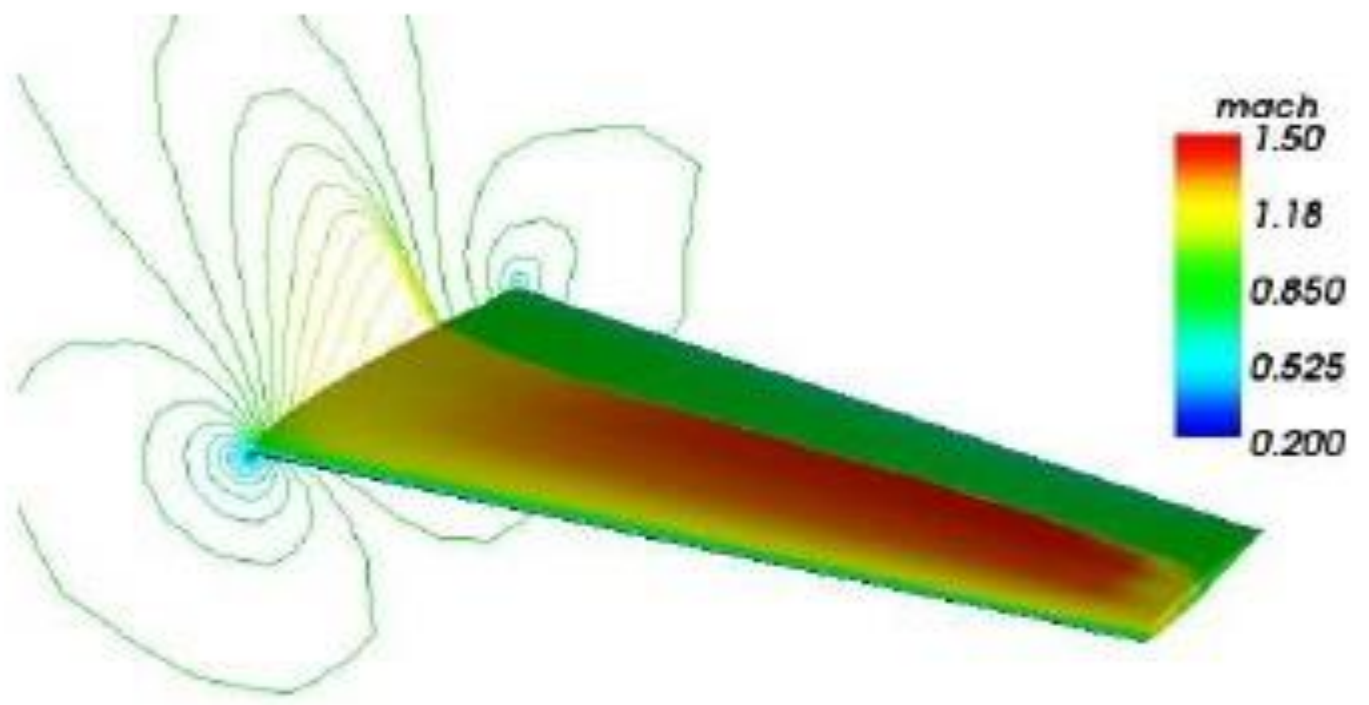

Figure 7: Mach number field on the wing and Mach number contours in the symmetry plane: initial shape.

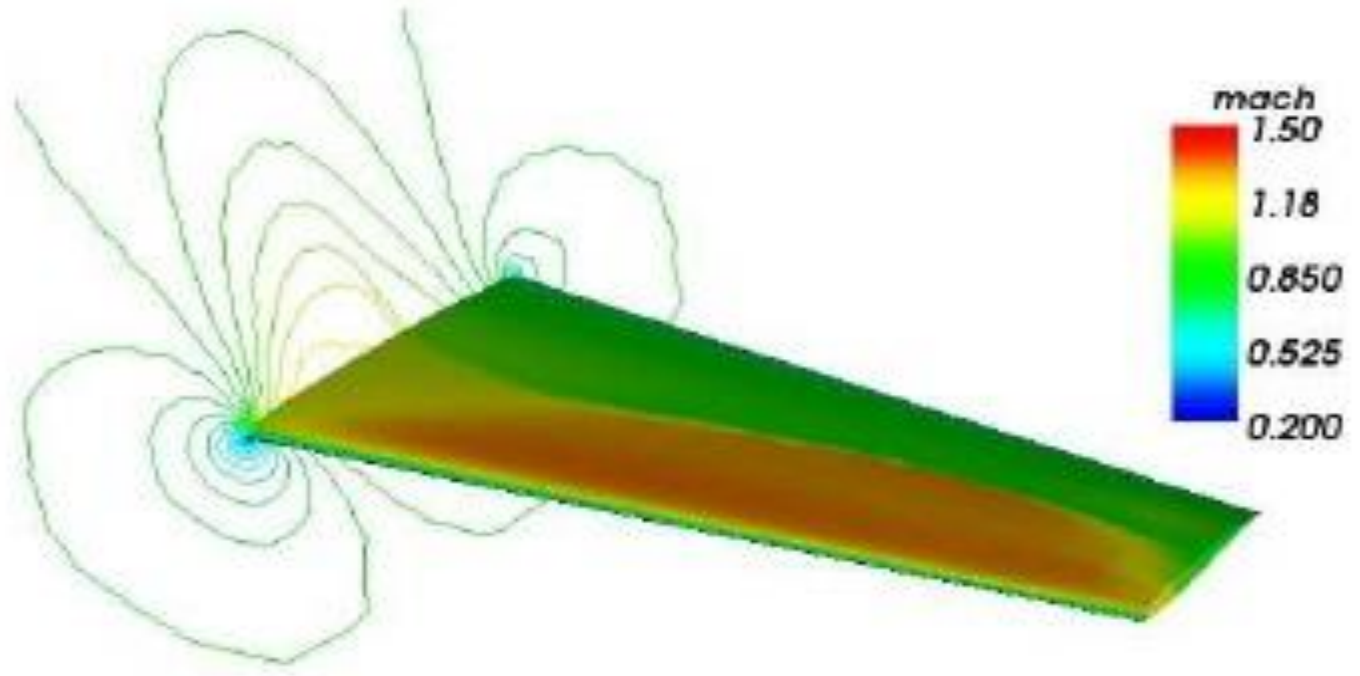

Figure 8: Mach number field on the wing and Mach number contours in the symmetry plane: basic method.

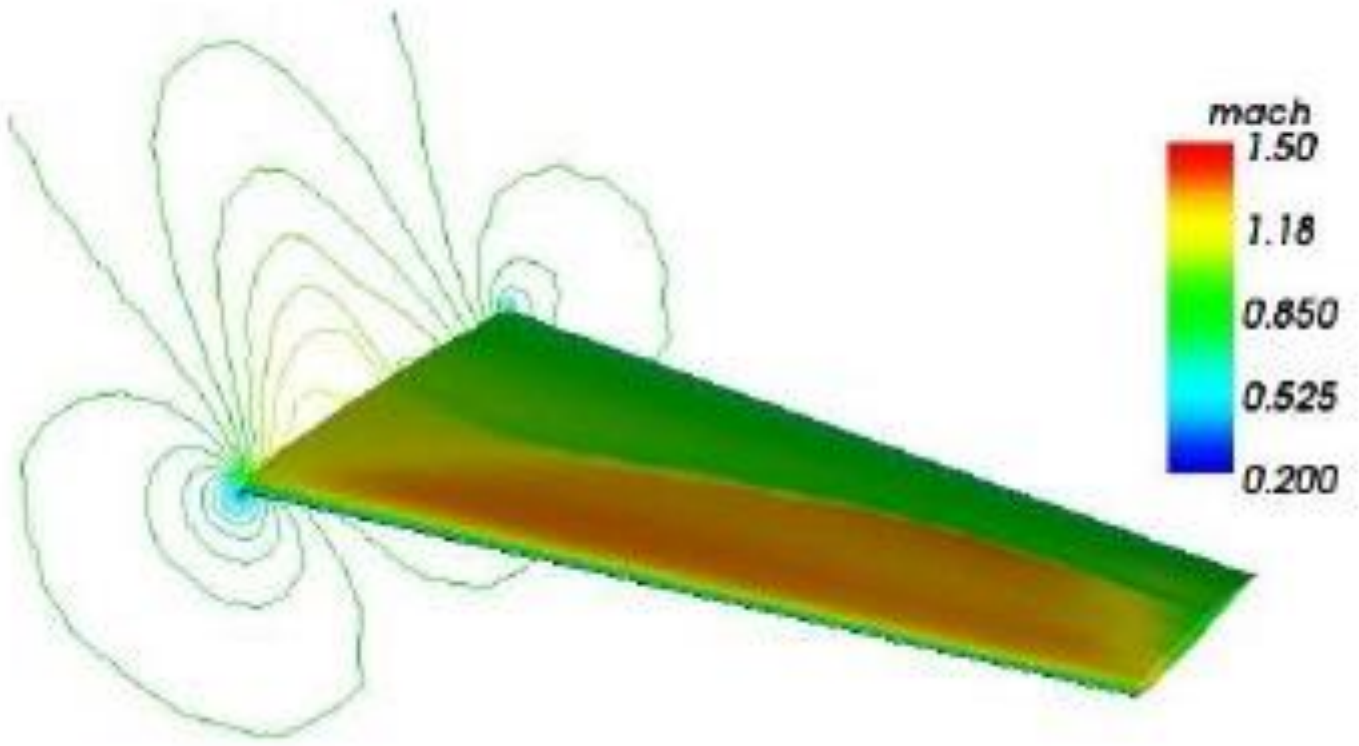

Figure 9: Mach number field on the wing and Mach number contours in the symmetry plane: adaptive method.

# 5. Conclusion

A self-adaptive procedure was developped for 3D optimum design in aerodynamics by using Free-Form Deformation (FFD). This approach which regularize the shape representation during the optimization process, is very effective in accelerating the convergence rate as shown by the numerical results. It follows clearly that the parametrization has a great impact on the results of aerodynamic shape optimization. this self-adaptive approach can be used to support a multilevel strategies developed in [4] [13][14].

This study can be extended for more complex engineering designs by using isogeometric analysis based on advanced techniques like non-uniform rational B-splines (NURBS) [2] due to its compact and shape representation capability [8]. Isogeometric analysis was introduced at first by Hughes et al. [9] to close the gap between computer aided design (CAD) and finite element analysis. Since then, many researchers in the fields of computational mechanical and geometric computation were involved in this topic [10][11][12]. The main idea behind isogeometric analysis is that the basis used to exactly parametrize the geometry will also serve as the basis for the solution space of the numerical method. This technique allows to best represent the shape and reduce the time required for its analysis. Although both positions and weights of control points affect the NURBS geometry instead of taking only the positions of control points as design variables when using the Bézier parameterization.

**Acknowledgements** The author wishes to express his warmest thanks to his collaborators in the INRIA Opale Project-Team, and particularly to J-A. Désidéri, A. Habbal and R. Duvigneau for their opinion and their valuable pieces of advice.

# References

[1] T.W. Sederberg and S.R. Parry, *Free-Form Deformation of solid geometric models*, Computer Graphics, 20(4), pp. 151-160, 1986. See also: Computer Aided geometric Design, T.W. Sederberg, ttp://tom.cs.byu.edu/"tom/ (Item: Courses).

[2] G. Farin, *Curves and Surfaces for Computer-Aided Geometric Design - A Practical Guide*, W. Rheinboldt and D. Siewiorek eds, Academic Press, Boston, 1990.

[3] J.-A. Désidéri, B. Abou El Majd, and A. Janka, *Nested self-adaptive Bézier parametrization for shape optimization*, Journal of Computational Physics, 224 (2007), pp. 117-131.

[4] B. Abou El Majd, R. Duvigneau, and J.-A. Désidéri, *Aerodynamic shape optimization using a full and adaptive multilevel algorithm*, ERCOFTAC 2006 Design Optimization : Methods & Application, 5-7 April 2006, Gran Canaria, Canaria Island, Spain.

[5] J.-A. Désidéri, *Two-level ideal algorithm for parametric shape optimization*, Advances in Numerical Mathematics, W. Fitzgibbon et al. eds, 2006.

[6] J. A. Samareh, *Multidisciplinary Aerodynamic-Structural Shape Optimization Using Deformation (MASSOUD)*, 8th AIAA/NASA/USAF/ISS, MO Symposium on Multidisciplinary Analysis and Optimization, September 6-8, 2000/Long Beach, CA, n! 4911 in AIAA-2000, 2000. J. A. Nelder, R. Mead, *A Simplex Method for Function Minimization*, Computer Journal, vol. 7, p. 308-313, 1965.

[7] P. C. Hansen, P. C. *Rank-Deficient and Discrete Ill-Posed Problems: Numerical Aspects of Linear Inversion*, SIAM (1998), Philadelphia.

[8] T.J.R. Hughes, J.A. Cottrell, Y. Bazilevs. *Isogeometric analysis: CAD, finite elements, NURBS, exact geometry, and mesh refinement*. Computer Methods in Applied Mechanics and Engineering, 194(2005) 4135-4195.

[9] G. Xu, B. Mourrain, R. Duvigneau, A. Galligo. *Optimal analysis-aware parameterization of computational domain in 3D isogeometric analysis*. Computer-Aided Design, 45(2013) 812-821. Y. Bazilevs, V.M. Calo, T.J.R. Hughes, and Y. Zhang. *Isogeometric fluid structure interaction: Theory, algorithms, and computations*. Computational Mechanics, 43(2008) 3-37.

[10] Y. Bazilevs, V.M. Calo, J.A. Cottrell, J. Evans, T.J.R. Hughes, S. Lipton, M.A. Scott, and T.W. Sederberg. *Isogeometric analysis using T-Splines*. Computer Methods in Applied Mechanics and Engineering, 199(2010) 229-263

[11] J. Zhao, J-A. Désidéri, and B. Abou El Majd. Two level correction algorithms for model problems, 2007.

[12] J. Zhao, B. Abou El Majd, and J-A. Désidéri. Two Level Correction Algorithm for Parametric Shape Inverse Optimization. International Journal of Engineering and Mathematical Modelling 2.1 (2015): 17-30.

[13] J-A. Désidéri, Jean-Antoine, R. Duvigneau, and B. Abou El Majd. Algorithms for efficient shape optimization in aerodynamics and coupled disciplines. 42nd AAAF Congress on Applied Aerodynamics, (Sophia-Antipolis, France), 2007.

[14] B. Abou El Majd, J-A. Désidéri, T. Do, L. Fourment, A. Habbal, and A. Janka. Multilevel strategies and hybrid methods for shape optimization and application to aerodynamics and metal forming. In : *Evolutionary and Deterministic Methods for Design, Optimization and Control with Applications to Industrial and Societal Problems Conference* (EUROGEN 2005), 2005.